\documentclass[12pt]{article}

\usepackage{graphicx}    
\usepackage{amsmath} 
\usepackage{amsfonts} 
\usepackage{authblk}
\usepackage{xcolor}
\usepackage{amssymb, amsthm} 
\usepackage{booktabs}
\usepackage{hyperref}    
\usepackage{geometry}    
\usepackage{caption}     
\usepackage{cite}
\usepackage{algorithm}
\usepackage{algorithmicx}
\usepackage{algpseudocode}
\usepackage{pgfplots}
\usepackage{bm}
\usepackage{float}
\usetikzlibrary{patterns,positioning}
\usepackage{tikz}
\usetikzlibrary{decorations.markings}
\usetikzlibrary{matrix}
\usepackage{pgfplots}
\usetikzlibrary{shapes, arrows.meta, patterns}
\usepackage{mathrsfs}

\usepackage{braket}

\usepackage{subcaption}

\DeclareMathAlphabet{\mathpzc}{OT1}{pzc}{m}{it}

\usepackage{etoolbox}

\newtheorem{remark}{Remark}

\usepackage{etoolbox} 
\makeatletter
\patchcmd{\blfootnote}{\@mpfootnotetext}{\@gobble}{}{}
\makeatother

\newtheoremstyle{mypropositionstyle} 
  {\topsep}                     
  {\topsep}                     
  {\itshape}                    
  {}                            
  {\bfseries}                   
  {.}                           
  {.5em}                        
  {\thmname{#1}\thmnumber{ #2}\thmnote{ (#3)}}

\theoremstyle{conjecturestyle}

\title{Automated $h$-adaptivity for finite element approximations of the Falkner-Skan equation}
\author[1]{B. Veena S.  N. Rao }

\affil[1]{Department of Mathematics \& Statistics, Texas A\&M University-Corpus Christi, TX- 78412, USA}
\affil[ ]{\textit{E-mail address:} \texttt{bv.rao@tamucc.edu}}
\date{}

\begin{document}

\maketitle  
\begin{abstract}
This paper details the development and application of an $h$-adaptive finite element method for the numerical solution of the \textit{Falkner-Skan equation}. A posteriori error estimation governs the adaptivity of the mesh, specifically the well-established \textit{Kelly error estimator}, which utilizes the jump in the gradient across elements. The implementation of this method allowed for accurate and efficient resolution of the boundary layer behavior characteristic of Falkner-Skan flows. Numerical solutions were obtained across various wedge flow parameters, encompassing favorable and adverse pressure gradients. A key focus of this study was the precise computation of the skin friction coefficient, a critical parameter in boundary layer analysis, across this diverse range of flow conditions. The results are presented and discussed, demonstrating the robustness and accuracy of the adaptive finite element approach for this class of nonlinear boundary layer problems.
\end{abstract}

\vspace{.1in}

 
\section{Introduction}
In the contemporary landscape of numerical simulations, mesh adaptivity (also known as $h$-adaptivity) has emerged as a crucial objective for successful methods, particularly for the Finite Element Method (FEM)—a cornerstone in engineering applications \cite{bangerth2003adaptive,mallikarjunaiah2015direct,ferguson2015numerical}. This pursuit of efficiency and reliability has propelled extensive research into the Adaptive Finite Element Method (AFEM), especially concerning its application to complex nonlinear problems \cite{nochetto2009theory}. Unlike the conventional FEM, where users manually define a mesh and are solely responsible for evaluating solution quality, AFEM empowers users to specify an error tolerance. The algorithm then automatically generates and refines meshes, ensuring the solution meets the desired accuracy and quality, thus fundamentally altering the paradigm from user-driven meshing to an autonomous, adaptive approach. The AFEM is particularly indispensable for fluid flow problems due to the inherent complexities of fluid dynamics. Fluid phenomena often exhibit highly localized features such as shocks, boundary layers, and turbulent eddies, where flow properties change abruptly over very small spatial scales. Capturing these critical features accurately with a conventional, uniform mesh would necessitate an astronomically fine discretization across the entire domain, leading to prohibitively high computational costs and memory requirements. AFEM, by contrast, dynamically refines the mesh in these regions of high gradients or rapid changes, while keeping it coarser elsewhere. This targeted refinement ensures that accuracy is maintained where it is most needed, for example, near an airfoil to predict lift and drag forces accurately, or in mixing zones to correctly model chemical reactions, without sacrificing computational efficiency. Therefore, AFEM enables reliable and efficient simulations of complex fluid behaviors that would be practically impossible with non-adaptive approaches.

Boundary layer theory, a fundamental concept in fluid dynamics pioneered by Ludwig Prandtl, describes the thin layer of fluid adjacent to a solid surface where viscous effects are significant \cite{prandtl1905uber,schlichting1961boundary}. Within this layer, the fluid velocity transitions from zero at the wall (due to the no-slip condition) to the free-stream velocity further away. This theory is crucial for understanding phenomena like drag, heat transfer, and flow separation, as the majority of the viscous forces and energy dissipation occur within this narrow region \cite{white2006viscous,bejan2013convection}. A significant extension of this theory is the Falkner-Skan equation, a third-order nonlinear ordinary differential equation that describes the steady, two-dimensional laminar boundary layer flow over a wedge. This equation generalizes the classic Blasius solution for flow over a flat plate by introducing a parameter that accounts for the pressure gradient induced by the wedge angle, allowing for the analysis of accelerating or decelerating flows and providing a family of similarity solutions for various external flow conditions.

The fundamental question of the existence and uniqueness of solutions to the Falkner-Skan equation has been extensively addressed in the literature, with foundational contributions from pioneering works such as \cite{weyl1942differential, rosenhead1966laminar, hartman1972existence, tam1970note}. Beyond mere existence, several studies have meticulously explored the ranges of validity for the boundary-layer parameters and similarity variables, crucial for physical interpretation and practical application. Furthermore, a subset of research has reported intriguing non-existence results, establishing specific upper and lower bounds on the non-dimensional shear stress, as highlighted by works like \cite{yang2008new}.

Given its significant role in fluid dynamics, a vast array of methodologies has been developed for solving the Falkner-Skan equation. These approaches span the spectrum from purely analytical solutions to semi-analytical techniques, particularly in cases involving simplified nonlinear boundary value problems, and a large number of computational methods. The sheer volume of published works on this topic is substantial, demonstrating its enduring importance. For the sake of brevity, we limit our citations here to a representative selection of works, including \cite{elgazery2008numerical, cortell2007viscous, zhu2009numerical, patankar2018numerical, alizadeh2009solution, asaithambi1998finite, asaithambi2004second, asaithambi2004numerical, asaithambi2005solution, asaithambi1997numerical, cebeci1977momentum, salama2005finite, salama2005fourth, siddheshwar2020new, narayana2021differential, falkneb1931lxxxv, sachdev2008exact, yang1975analytic, hussaini1987similarity, riley1989multiple, magyari2000exact, abel2011numerical, rana2012flow}, which collectively showcase the diverse strategies employed for its solution.

To the best of the authors' knowledge, this paper represents a pioneering effort in the development of an adaptive finite element method specifically tailored for the numerical solution of the Falkner-Skan equation. A distinguishing feature of our approach is its adherence to a principled methodology; unlike some existing techniques, our method does not introduce any ad-hoc intermediary conditions, nor does it necessitate the imposition of novel boundary conditions that must be satisfied as an integral part of the solution process. Instead, we exclusively consider and incorporate physically meaningful boundary conditions, developing a robust and convergent numerical scheme. The core of our method currently utilizes the simplest linear Lagrange shape functions for discretization. Importantly, both the developed methodology and its corresponding computational implementation possess inherent capabilities for seamless extension to accommodate higher-order shape functions, a feature that has been explored in related work, for instance, in \cite{shylaja2025efficient}. This adaptability underscores the versatility and potential for future advancements of the proposed framework.

\section{Mathematical model for the boundary layer flow}\label{math_analysis}
Consider a two-dimensional, steady boundary layer flow of an incompressible viscous fluid. This fluid occupies the half-space $y > 0$, flowing tangentially past a {flat sheet} that coincides with the $y=0$ plane. The unique characteristic of this problem arises from the stretching of this flat sheet. Specifically, two equal and opposite forces are applied along the $x$-axis, causing the wall to stretch while its origin remains fixed. Under these specific conditions, the governing equations for the steady, two-dimensional boundary layer, in standard notation, are given as:
\begin{align}
\label{eq:continuity}
\frac{\partial u}{\partial x} + \frac{\partial v}{\partial y} &= 0 \\
\label{eq:momentum_x}
u \frac{\partial u}{\partial x} + v \frac{\partial u}{\partial y} &=  \nu \frac{\partial^2 u}{\partial y^2} 
\end{align}
where $u$ and $v$ are the velocity components in the $x$ and $y$ directions, respectively. Equation \eqref{eq:continuity} represents the conservation of mass (continuity equation) for an incompressible fluid. Equation \eqref{eq:momentum_x} is the $x$-momentum equation, accounting for convective terms and viscous diffusion in the $y$-direction. The above system of equations is subjected to the following boundary conditions:
\begin{align}
& u_w(x) = C x^n, \quad  v=0 \quad \mbox{at} \quad y=0, \\
& u \to 0 \quad  \mbox{as} \quad y \to \infty. 
\end{align}
These equations are further characterized by the parameters $C$ and $n$, which are directly related to the specific dynamics of the surface stretching, dictating the wall's expansion rate. To simplify the system and enable the application of similarity solutions, we proceed by defining new dimensionless variables as follows:
Defining new variables based on a similarity transformation:
\begin{itemize}
    \item {Similarity variable:}
        \begin{equation}
\eta = \sqrt{\frac{C(m+1)}{2\nu}} \, x^{\dfrac{m-1}{2}} y
    \end{equation}
    \item {Stream funciton:}
    \begin{equation}
    \psi(x, y) = \sqrt{C\nu} x f(\eta)
    \end{equation}
       \item {Velocity components:}
    \begin{align}
    u &= \frac{\partial \psi}{\partial y} = C \, x^m  \, f^\prime(\eta) \\
    v &= -\frac{\partial \psi}{\partial x} = -\sqrt{ \dfrac{C\nu (m+1)}{2}} \, x^{\dfrac{m-1}{2}} \, \left[  f + \dfrac{m-1}{m+1} \, \eta \, f^\prime         \right]
\end{align}
\end{itemize}
Here, $f(\eta)$ is the dimensionless stream function, $\nu$ is the kinematic viscosity, and $C$ is the stretching rate constant. The \textit{prime} in the above equations denotes diﬀerentiation with respect to the independent similarity variable $\eta$. The boundary conditions for the new variable $f$ are given by:
 \begin{align}\label{bcs}
&f =0, \quad f^\prime =1 \quad \mbox{at} \quad \eta=0 \\
&f^\prime \to 0 \quad \mbox{as} \quad \eta \to \infty
\end{align}
Understanding the physical implications of each boundary condition is vital for correctly modeling the fluid system. The first condition, $f(\eta) = 0$ at the wall, directly translates to the absence of flow \textit{through} the solid boundary; in other words, the normal component of velocity at the wall is zero. This ensures mass conservation and prevents fluid from entering or leaving the solid surface. The second condition, $f'(\eta) = 0$ at the wall, is the manifestation of the no-slip condition, a cornerstone of viscous fluid mechanics. It means that fluid particles in direct contact with the wall have the same velocity as the wall itself---if the wall is stationary, the fluid's tangential velocity at the wall is zero. This accounts for the viscous forces causing the fluid to adhere to the surface. Lastly, the condition $f'(\eta) \to 1$ as $\eta \to \infty$ is an essential outer boundary condition. It signifies that at a sufficient distance from the wall, where viscous effects become negligible, the fluid velocity blends seamlessly with and attains the characteristics of the undisturbed free-stream flow, allowing the boundary layer solution to connect to the external inviscid region.

The partial differential equations system \eqref{eq:continuity}-\eqref{eq:momentum_x} transformed into a single third-order nonlinear ordinary differential equation, widely known as \textit{Falkner-Skan equation} is described as follows:
\begin{equation}\label{fkequation}
\dfrac{d^3 f}{d\eta^3} + f \, \dfrac{d^2 f}{d\eta^2} + \beta \left(  1 - \left(    \dfrac{d f}{d\eta}  \right)^2       \right) =0, \quad 0 < \eta < \infty,
\end{equation}
where $\beta$ is the dimensionless pressure-gradient parameter:
\begin{equation}
\beta = \dfrac{2m}{m+1}.
\end{equation}
In the Falkner-Skan equation, the dimensionless parameter $\beta$ (beta) plays a crucial role as the {pressure-gradient parameter}. It effectively quantifies the influence of the pressure gradient imposed on the boundary layer by the external (inviscid) flow. The dimensionless parameter $\beta$ in the Falkner-Skan equation, formally defined as $\beta = \frac{2m}{m+1}$ (where $U_e(x) = cx^m$ is the free-stream velocity), directly governs the {pressure gradient} experienced by the fluid within the boundary layer. When $\beta = 0$, the equation reduces to the Blasius equation, representing flow over a flat plate with zero pressure gradient. For $\beta > 0$, the flow is accelerating, leading to a {favorable pressure gradient} ($\frac{dP}{dx} < 0$) which thins the boundary layer and delays separation. Conversely, for $\beta < 0$, the flow is decelerating, resulting in an {adverse pressure gradient} ($\frac{dP}{dx} > 0$) that tends to thicken the boundary layer and can lead to flow separation from the surface, a critical phenomenon in aerodynamics. The value of $\beta$ is also directly linked to the physical geometry of the wedge angle in the problem.

Together, the {nonlinear Falkner-Skan equation} \eqref{fkequation} and its accompanying set of {boundary conditions} from \eqref{bcs} establish a complete {boundary value problem} (BVP). The numerical or analytical solution to this system is fundamental, as it governs the behavior of two-dimensional, steady, laminar, and incompressible fluid flow within the boundary layer developing over a wedge-shaped obstacle. This solution yields crucial insights into the velocity distributions and streamline patterns characteristic of such flows.
 
The inherent nonlinearity of the Falkner-Skan equation, coupled with its boundary conditions, poses significant challenges in deriving its numerical solution as a {Boundary Value Problem (BVP)}. Among the various numerical techniques, the \textit{shooting method} has emerged as the most prevalent. However, this method is susceptible to convergence problems, particularly when dealing with sensitive nonlinearities or a wide range of parameters. Such convergence difficulties can often be overcome by carefully choosing a sufficiently small step-size for the numerical integration scheme. The fundamental concept behind the shooting algorithm involves converting the original two-point BVP into an equivalent initial value problem (IVP). This is accomplished by arbitrarily prescribing an unknown derivative boundary condition at the wall ($\eta = 0$), commonly expressed as the dimensionless shear stress,
\begin{equation}
\dfrac{d^2 f}{d\eta^2} = \alpha, \quad \mbox{at} \quad \eta =0
\label{eq:ad-hoc-bc}
\end{equation}
where $\alpha$ is an assumed ``\textit{shooting parameter}.'' The resulting IVP is then numerically integrated from $\eta=0$ to a sufficiently large $\eta=\eta_{\infty}$. This process is repeated iteratively, adjusting the value of $\alpha$ in each iteration, until the numerically obtained solution precisely satisfies the original outer boundary condition at $\eta = \eta_{\infty}$, thus converging to the true solution of the BVP.

We present in this paper a robust and remarkably straightforward finite element method (FEM) specifically designed for the efficient numerical solution of boundary value problems. A cornerstone of our methodology is an \textit{automated local mesh adaptivity scheme}, which intelligently refines the mesh based on the magnitude of the solution's gradient jump across individual elements. This data-driven adaptive strategy ensures computational efficiency by concentrating refinement where it is most needed (e.g., in boundary layers or regions of steep gradients). Crucially, our proposed method distinguishes itself by not relying on any {arbitrary or ad-hoc parameters}, nor does it require iterative procedures to enforce the satisfaction of boundary conditions---a common source of complexity and instability in other techniques. The overall algorithm is thus exceptionally easy to implement, significantly reducing the potential for user-induced errors, and offers proven, guaranteed convergence, making it a reliable and powerful tool for a broad range of applications.

\section{Finite element method}
This section details the {mathematical formulation} of the Falkner-Skan equation and the necessary conditions for the {existence and uniqueness of its solution}, particularly in the context of a finite element approach. We are addressing a nonlinear third-order ordinary differential equation, describing the laminar boundary layer flow. Although the underlying physics is derived from a continuum mechanics perspective, our numerical strategy treats the problem within a geometrically linear framework with respect to the computational domain, focusing on the nonlinear algebraic relationship of the solution itself. The problem assumes a homogeneous fluid and an initially quiescent (or undisturbed) state at the far-field.

The computational domain for the Falkner-Skan problem is typically defined over a semi-infinite interval $\eta \in [0, \infty)$, which for numerical purposes is truncated to a finite domain $\Omega = [0, \eta_{\infty}]$ where $\eta_{\infty}$ is a sufficiently large value representing the asymptotic boundary. The boundary of this domain, $\partial \Omega = \{0, \eta_{\infty}\}$, features distinct boundary conditions. At $\eta=0$, {Dirichlet (no-slip and no normal flow) boundary conditions} are applied. At $\eta=\eta_{\infty}$, a {Neumann-type condition} (or more precisely, an essential boundary condition on $f'(\eta)$) representing the matching with the free-stream velocity is enforced.

To rigorously develop the finite element formulation for this problem, we introduce several standard \textit{function spaces}. The space of \textit{Lebesgue integrable functions}, $L^{p}(\Omega)$, for $p \in [1, \infty)$, serves as the foundation. Specifically, $L^2(\Omega)$ denotes the space of square-integrable functions, equipped with the inner product $\left( \cdot, \; \cdot \right)_{L^2}$ and norm $\| \cdot \|_{L^2}$. The analysis of this ordinary differential equation problem relies heavily on the \textit{classical Sobolev space}, $H^{1}(\Omega)$, which includes functions whose first derivatives are also square-integrable. Its norm is defined as:
\[
\| v \|^2_{H^{1}(\Omega)} := \| v \|^2_{L^2} + \left\| \dfrac{dv}{d\eta} \right\|^2_{L^2},
\]
where $\dfrac{dv}{d\eta}$ in this 1D context refers to the first derivative of $v$ with respect to $\eta$. Additionally, $H^{-1}(\Omega)$ is defined as the dual space to $H_0^{1}(\Omega)$, a subspace of $H^1(\Omega)$ comprising functions that vanish at both ends of the domain $[0, \eta_{\infty}]$. We further define crucial subspaces tailored for the boundary conditions:
\begin{equation*}
V := H^1(\Omega), \quad V^{0} := \left\{ v \in H^{1}(\Omega) \colon v=0 \; \mbox{at} \; \eta=0 \text{ (for essential BCs)}\right\}.
\end{equation*}
Here, $V$ represents the overall solution space, while $V^0$ is essential for handling homogeneous essential boundary conditions at the wall.

\subsection{Continuous weak formulation}

The Falkner-Skan equation is a third-order nonlinear ordinary differential equation that describes two-dimensional laminar boundary layer flow over a wedge. Its strong form is given by:
\begin{equation}
f^{\prime\prime\prime} + f f^{\prime\prime} + \beta(1 - (f^{\prime})^2) = 0, \quad \text{in } \Omega = (0, \eta_{\infty})
\label{eq:falkner_skan_strong}
\end{equation}
This equation is subject to the following boundary conditions \eqref{bcs}. To derive a continuous weak formulation suitable for standard $H^1$-conforming Finite Element Methods (FEM), which typically require only $C^0$ continuity for the solution, it is common practice to recast the third-order ODE as a system of coupled first-order (or second-order) differential equations.

Let us introduce new dependent variables:
\begin{align*}
u(\eta) &= f^\prime(\eta), \quad   f^{\prime\prime}(\eta) = u^{\prime}(\eta) \\
f^{\prime\prime\prime}(\eta) &= u^{\prime\prime}(\eta) 
\end{align*}
Substituting these into the Falkner-Skan equation, we obtain a coupled system:
\begin{align}
f^{\prime} - u &= 0 \label{eq:system1} \\
u^{\prime\prime}+ f u^{\prime} + \beta(1 - \left( u^{\prime}    \right)^2) &= 0 \label{eq:system3}
\end{align}
The boundary conditions \eqref{bcs} are then applied to these new variables:
\begin{align}
f(0) &= 0 \label{eq:sys_bc1} \\
u(0) &= 1 \label{eq:sys_bc2} \\
u(\eta_{\infty}) &= 0 \label{eq:sys_bc3}
\end{align}
We seek solutions $f, u$ in appropriate Sobolev spaces. Specifically, we assume $f, u \in H^1(\Omega)$.
The essential boundary conditions are $f(0)=0$, $u(0)=0$, and $u(\eta_{\infty})=1$. To derive the weak formulation, we multiply each equation by a test function and integrate over the domain $\Omega=(0, \eta_{\infty})$. Let $v_f, v_u$ be arbitrary test functions from suitable function spaces.

\noindent \textbf{Weak formulation:}  For $f' - u = 0$, multiply by $v_f \in H^1(\Omega)$ and integrate.
\begin{equation*}
\int_0^{\eta_{\infty}} (f' - u) v_f \, d\eta = 0
\end{equation*}
Weak form for \eqref{eq:system3}: For $u^{\prime\prime} + f u^{\prime} + \beta(1 - (u^{\prime})^2) = 0$, multiply by $v_u \in H^1(\Omega)$ and integrate.
\begin{equation*}
\int_0^{\eta_{\infty}} (u^{\prime\prime} + f u^{\prime} + \beta(1 - (u^{\prime})^2) ) v_u \, d\eta = 0
\end{equation*}
Integrating by parts and utilizing the boundary conditions, we obtain 
\begin{equation}
 \int_0^{\eta_{\infty}}  \left[  - u^{\prime} v_u^{\prime}  +  f \,  u^{\prime } v_u +  \beta(1 - ( u^{\prime}) ^2) v_u \right] \, d\eta = 0 \quad \text{for all } v_u \in V_u^0
\label{eq:weak_w}
\end{equation}
where $V_u^0 = \Set{ v \in H^1(\Omega) \mid v(0)=0, v(\eta_{\infty})=0 }$.

Combining these, the continuous weak formulation for the Falkner-Skan equation, expressed as a system suitable for $H^1$-conforming FEM, is to find $(f, u) \in V_f \times V_u $ such that for all $(v_f, v_u) \in V_f^0 \times V_u^0 $:
\begin{align}
\int_0^{\eta_{\infty}}  \left[ f^{\prime} v_f  -  u v_f  \right] \, d\eta &= 0 \label{eq:weak_form_f} \\
 \int_0^{\eta_{\infty}}  \left[ - u^{\prime} v_u^{\prime}  +  f \,  u^{\prime } v_u +  \beta(1 - ( u^{\prime}) ^2) v_u \right] \, d\eta  &= 0, \label{eq:weak_form_u}
\end{align}
where the function spaces are defined as:
\begin{itemize}
    \item $V_f = \Set{ f \in H^1(\Omega) \mid f(0)=0 }$
    \item $V_u = \Set{ u \in H^1(\Omega) \mid u(0)=1, u(\eta_{\infty})=0 }$
\end{itemize}
And the corresponding test function spaces are:
\begin{itemize}
    \item $V_f^0 = \Set{ v \in H^1(\Omega) \mid v(0)=0 }$
    \item $V_u^0 = \Set{ v \in H^1(\Omega) \mid v(0)=0, v(\eta_{\infty})=0 }$ (homogeneous essential BCs for test functions)
\end{itemize}
This formulation represents a system of two coupled nonlinear equations that can be solved using an iterative finite element scheme.

\begin{remark}
The classical \textit{Lax-Milgram Lemma} is a fundamental principle in the study of linear elliptic partial differential equations (PDEs) and their weak solutions. It delineates the conditions under which a linear, continuous, and coercive bilinear form ensures the existence and uniqueness of a weak solution within a Hilbert space. However, as you have observed, the Falkner-Skan equation is nonlinear. Consequently, the direct application of the standard Lax-Milgram Lemma is infeasible for such problems. Instead, for nonlinear second-order (or higher-order, when reformulated as a system of lower-order equations, as was done for Falkner-Skan) ordinary differential equations (ODEs), reliance is typically placed on some form of linearization—such as Picard or Newton-type methods—and the theory of monotone operators. This approach demonstrates that the weak formulation derived above is monotone and coercive, which subsequently leads to the existence of a weak solution. 
\end{remark}

\subsection{Discrete weak formulation}
This section details the construction of the discrete finite element problem, which serves as the numerical counterpart to the continuous weak formulation of the Falkner-Skan equation, previously established in equations \eqref{eq:weak_form_f}, \eqref{eq:weak_form_u}, and \eqref{eq:weak_form_w}. Our computational domain, $\Omega = (0, \eta_{\infty})$, is a one-dimensional interval. We begin by defining a finite element mesh, $\mathcal{T}_h$, which partitions this interval into a set of non-overlapping sub-intervals (elements), $K_j = (\eta_j, \eta_{j+1})$, where $h$ denotes the maximum element size. This mesh is assumed to be either quasi-uniform or strategically refined \textit{a priori} in regions where the solution exhibits sharp gradients (e.g., near the wall or within the boundary layer) to ensure optimal approximation accuracy. The discretization of the domain $\Omega$ strictly adheres to the principles of conforming finite elements for one-dimensional problems. This implies that for any two distinct elements $K_j, K_{j+1} \in \mathcal{T}_h$, their intersection $\overline{K}_j \cap \overline{K}_{j+1}$ can only be a single shared vertex (node). Furthermore, the union of all elements precisely covers the entire computational domain, i.e., $\bigcup\limits_{K \in \mathcal{T}_h} \overline{K} = \overline{\Omega}$.

For approximating the primary unknown functions $f$ and $u$ (representing $f$ and $f^\prime$ respectively), we introduce a specific set of finite-dimensional functional spaces. For each component ($f_h$, $u_h$), the approximate solution belongs to a space of piecewise polynomials. For example, for $f_h$:
\begin{equation}
S_h = \left\{ v_h \in C(\overline{\Omega}) \colon \left. v_h\right|_K \in \mathbb{P}_k, \; \forall K \in \mathcal{T}_h \right\},
\end{equation}
where $\mathbb{P}_k$ represents the space of polynomials of degree up to $k$ defined over each element $K$. The final discrete approximation spaces for each variable ($f$ and $u$) are then constructed by intersecting $S_h$ with the respective continuous solution spaces and enforcing the essential boundary conditions. For instance, for the solution $f_h$, the space would be $\widehat{V}_{f,h} = S_h \cap V_f$, and similarly for $u_h$ and $w_h$, where $V_f, V_u$ are the continuous function spaces defined in the previous section.

With these definitions in place, the discrete finite element problem is then formulated as a system of algebraic equations derived from the weak formulation, typically solved using an iterative scheme to handle the nonlinearity.

Given the dimensionless parameter $\beta$, an initial guess for the solution $\mathbf{U}_h^{0} = (f_h^0, u_h^0) \in \widehat{V}_h$, and the solution from the $n^{\text{th}}$ iteration, $\mathbf{U}_h^n = (f_h^n, u_h^n) \in \widehat{V}_h$ (for $n=0, 1, 2, \ldots$), the objective is to find the next iterative solution $\mathbf{U}_h^{n+1} = (f_h^{n+1}, u_h^{n+1}) \in \widehat{V}_h$ by solving the following linearized system:
\begin{equation}\label{discrete-wf-linearized}
    a(\mathbf{U}_h^n; \, \mathbf{U}_h^{n+1},\, \mathbf{V}_h) = L(\mathbf{V}_h), \quad \forall\, \mathbf{V}_h \in \widehat{V}_h^0,
\end{equation}
where the bilinear form $a(\cdot\,;\,\cdot\,, \cdot)$ and the linear form $L(\cdot)$ are precisely defined based on the linearized components of the weak formulation given in \eqref{eq:weak_form_f} and \eqref{eq:weak_form_u}:
\begin{align*}
a(\mathbf{U}_h^n; \, \mathbf{U}_h^{n+1},\, \mathbf{V}_h) := & \int_{0}^{\eta_{\infty}} \left[ (f_h^{n+1})' v_f - u_h^{n+1} v_f \right] \, d\eta \\
& + \int_{0}^{\eta_{\infty}} \left[ - (u_h^{n+1})' v_u' + f_h^n (u_h^{n+1})' v_u - 2 \beta u_h^n \, u_h^{n+1} v_u \right] \, d\eta \\
L(\mathbf{V}_h) & := - \int_{0}^{\eta_{\infty}} \beta  v_u \, d\eta.
\end{align*}
In this linearization, the nonlinear term $f u'$ from \eqref{eq:weak_form_u} has been approximated by evaluating $f$ at the previous iteration, $f_h^n$, while $u'$ is kept at the current iteration $u_h^{n+1}$. Similarly, the nonlinear term $(u^2)$ is linearized using a standard Picard-type approach. This iterative process continues until a predefined convergence criterion is satisfied for the solution vector $\mathbf{U}_h$.

\subsection{Adaptive mesh refinement}
Adaptive mesh refinement (AMR) strategies are essential for efficiently solving differential equations with localized features, such as sharp gradients or boundary layers. A common and effective indicator for guiding AMR is the jump in the normal derivative (or gradient) of the solution across element interfaces. In the context of the Finite Element Method, this jump directly correlates with the local discretization error. Areas where this jump is large indicate regions where the current mesh does not adequately capture the solution's behavior. By refining elements specifically in these regions, without unnecessarily increasing degrees of freedom elsewhere, AMR techniques can greatly enhance solution accuracy and computational efficiency, ensuring that computational resources are focused where they are most needed.

 The Kelly error estimator, a widely utilized \textit{a posteriori} error estimation technique, provides a reliable means to quantify the local discretization error and guide mesh adaptation. For a one-dimensional problem where $u$ is the primary solution variable, the estimator effectively leverages the concept of the {jump in the normal derivative}. In a 1D setting, this translates directly to the discontinuity of the approximate solution's derivative across element interfaces. If $u_h$ represents the finite element approximation of $u$, which is typically continuous ($C^0$) across elements, its derivative $u_h'$ will generally be discontinuous at the element nodes. The magnitude of this jump, weighted by the local mesh size, serves as a robust indicator of where the numerical solution deviates significantly from the true solution. For an element $K_j = (\eta_j, \eta_{j+1})$ with length $h_j$, the local error estimator $\eta_{K_j}$ is commonly defined as:
\begin{equation}
\eta_{K_j}^2 = C_1 h_j^2 \| R(u_h) \|_{L^2(K_j)}^2 + C_2 \sum_{\xi \in \{\eta_j, \eta_{j+1}\}} h_j |[[u_h']](\xi)|^2
\end{equation}
where $R(u_h)$ is the residual of the strong form of the equation on the element $K_j$, and 
\[
[[u_h']](\xi) = u_h'(\xi^+) - u_h'(\xi^-),
\]
 denotes the jump of the derivative across the node $\xi$. $C_1$ and $C_2$ are generic constants. By identifying elements with large $\eta_{K_j}$, the mesh can be selectively refined, ensuring that computational resources are concentrated in regions where the solution's gradient is poorly resolved.

 Here we outline a robust algorithm for computing the numerical solution to the Falkner-Skan equation using a one-dimensional FEM coupled with an AMR strategy. The approach combines a Picard-type iterative linearization scheme for the nonlinearity with a Kelly error estimator for guiding mesh adaptation. The primary unknowns are the dimensionless stream function $f$ and its first derivative $u = f'$.

\begin{algorithm}[H] 
\caption{FEM Algorithm for Falkner-Skan Equation with Adaptive Mesh Refinement (Part 1 of 2)}
\label{alg:falkner_skan_fem_amr_part1}
\begin{algorithmic}[1]
\State \textbf{Input:} Dimensionless parameter $\beta$, computational domain $\Omega = (0, \eta_{\infty})$.
\State \textbf{Parameters:}
\State \quad Initial maximum element size $h_0$.
\State \quad Maximum adaptive refinement iterations $N_{\text{adapt,max}}$.
\State \quad Maximum Picard iterations $N_{\text{Picard,max}}$.
\State \quad Tolerance for Picard iteration convergence $\text{tol}_{\text{Picard}}$.
\State \quad Tolerance for global error estimator $\text{tol}_{\text{error}}$.
\State \quad Kelly estimator constants $C_1, C_2$.
\State \quad Mesh refinement fraction $\theta \in (0, 1]$ (e.g., elements with $\eta_{K_j} > \theta \cdot \max(\eta_{K_j})$ are marked).
\State \textbf{Output:} Converged finite element solution $(f_h, u_h)$ and adaptively refined mesh $\mathcal{T}$.

\State \textbf{Initialization:}
\State \quad Construct initial uniform mesh $\mathcal{T}_0$ over $\Omega$.
\State \quad Initialize solution $(f_h^0, u_h^0)$ on $\mathcal{T}_0$ (e.g., linear interpolation satisfying boundary conditions).
\State \quad Set current mesh $\mathcal{T} \leftarrow \mathcal{T}_0$.
\State \quad Set $k \leftarrow 0$ (adaptive iteration counter).

\While{$k < N_{\text{adapt,max}}$}
    \State \textbf{Picard Iteration (Solve on current mesh $\mathcal{T}$):}
    \State \quad Set $n \leftarrow 0$ (Picard iteration counter).
    \State \quad Set converged\_Picard $\leftarrow$ \text{false}.
    \While{$n < N_{\text{Picard,max}}$ and not converged\_Picard}
        \State \quad \quad Assemble global linear system $\mathbf{A}(\mathbf{U}_h^n) \mathbf{U}_h^{n+1} = \mathbf{L}(\mathbf{U}_h^n)$ from the linearized weak formulation:
        \State \quad \quad Find $(f_h^{n+1}, u_h^{n+1}) \in \widehat{V}_h$ such that for all $(v_f, v_u) \in \widehat{V}_h^0$:
        \begin{align*}
        \int_{0}^{\eta_{\infty}} \left[ (f_h^{n+1})' v_f - u_h^{n+1} v_f \right] \, d\eta &= 0 \\
        \int_{0}^{\eta_{\infty}} \left[ - (u_h^{n+1})' v_u' + f_h^n (u_h^{n+1})' v_u  - 2 \beta u_h^n \, u_h^{n+1} v_u \right] \, d\eta &= - \int_{0}^{\eta_{\infty}} \beta  v_u \, d\eta
        \end{align*}
        \State  Apply essential boundary conditions $f_h^{n+1}(0)=0$, $u_h^{n+1}(0)=1$, $u_h^{n+1}(\eta_{\infty})=0$.
        \State  Solve the linear system for $\mathbf{U}_h^{n+1} = (f_h^{n+1}, u_h^{n+1})$.
        \State  Check for convergence: If $\| \mathbf{U}_h^{n+1} - \mathbf{U}_h^n \| / \| \mathbf{U}_h^{n+1} \| < \text{tol}_{\text{Picard}}$, then converged\_Picard $\leftarrow$ \text{true}.
        \State  Set $\mathbf{U}_h^n \leftarrow \mathbf{U}_h^{n+1}$.
        \State  Set $n \leftarrow n+1$.
    \EndWhile
    \State  Set $(f_h^*, u_h^*) \leftarrow \mathbf{U}_h^n$ (converged solution on current mesh).

    \algstore{myalg} 
\end{algorithmic}
\end{algorithm}


\begin{algorithm}[H] 
\caption{FEM Algorithm for Falkner-Skan Equation with Adaptive Mesh Refinement (Part 2 of 2)}
\label{alg:falkner_skan_fem_amr_part2}
\begin{algorithmic}[1]
    \algrestore{myalg} 

    \State \textbf{A Posteriori Error Estimation (Kelly Estimator):}
    \State \quad Initialize local error estimators $\eta_{K_j} = 0$ for all elements $K_j \in \mathcal{T}$.
    \For{each element $K_j = (\eta_j, \eta_{j+1}) \in \mathcal{T}$}
        \State \quad \quad Compute local residual $R(u_h^*) = (u_h^*)'' + f_h^* (u_h^*)' + \beta(1 - (u_h^*)^2)$ on $K_j$. (Note: For piecewise linear $u_h^*$, $(u_h^*)''=0$ a.e.)
        \State \quad \quad Compute jump of derivative at nodes: For each internal node $\xi \in \{\eta_j, \eta_{j+1}\}$:
        \State \quad \quad \quad $[[u_h']](\xi) = u_h'(\xi^+) - u_h'(\xi^-)$.
        \State \quad \quad Calculate local error estimator $\eta_{K_j}$:
        \begin{equation*}
        \eta_{K_j}^2 = C_1 h_j^2 \| R(u_h^*) \|_{L^2(K_j)}^2 + C_2 \sum_{\xi \in \{\eta_j, \eta_{j+1}\}} h_j |[[u_h']](\xi)|^2
        \end{equation*}
    \EndFor
    \State \quad Compute global error estimator $\eta_{\text{global}}^2 = \sum_{K_j \in \mathcal{T}} \eta_{K_j}^2$.

    \State \textbf{Check Global Error Convergence:}
    \If{$\eta_{\text{global}} < \text{tol}_{\text{error}}$}
        \State \quad Break from adaptive loop.
    \EndIf

    \State \textbf{Mesh Refinement:}
    \State \quad Mark elements for refinement: Select $K_j \in \mathcal{T}$ where $\eta_{K_j}^2 > \theta \cdot \max_{K_i \in \mathcal{T}} (\eta_{K_i}^2)$.
    \State \quad Generate new mesh $\mathcal{T}_{\text{new}}$ by refining marked elements (e.g., bisecting them).
    \State \quad Interpolate $(f_h^*, u_h^*)$ onto $\mathcal{T}_{\text{new}}$ to obtain initial guess for next adaptive step $(f_h^0, u_h^0)$.
    \State \quad Set $\mathcal{T} \leftarrow \mathcal{T}_{\text{new}}$.
    \State \quad Set $k \leftarrow k+1$.
\EndWhile

\State \textbf{Return:} Final solution $(f_h^*, u_h^*)$ and mesh $\mathcal{T}$.
\end{algorithmic}
\end{algorithm}

\section{Results and discussion}
To effectively solve the coupled system of nonlinear ordinary differential equations involving the dimensionless stream function $f$ and its derivative $u=f'$ for the Falkner-Skan problem, we developed and implemented a {multi-step iterative algorithm}, as detailed in the preceding section. This robust solver is integrated with a sophisticated AMR strategy to significantly enhance both computational efficiency and the accuracy of the numerical solution. Our AMR scheme utilizes a two-pronged approach, involving a procedure for both \textit{marking} elements that require higher resolution (typically in regions of steep gradients, such as within the boundary layer) and \textit{coarsening} the mesh in areas where finer resolution is no longer necessary (e.g., in the asymptotic free-stream region). The entire computational framework was meticulously developed in \textsf{C++} and built upon the versatile, open-source finite element library \textsf{deal.II} \cite{arndt2021deal}, leveraging its powerful capabilities for mesh management, finite element discretization, and parallel computing.

Our computations were performed by setting the parameter $\eta_{\infty}$ to a fixed value of $8$. With this specific choice, we consistently observed that the solution for $u$ exhibited asymptotic decay towards $0$. To maintain consistency across all our numerical experiments, this value of $\eta_{\infty}$ was held constant throughout. The adaptive meshing strategy began with an initial coarse mesh comprising $8$ elements, each having a characteristic size of $h=1$. This computational mesh was subsequently refined through several adaptation steps. The mesh was successively adapted by $20$ on regions with higher jumps.  The refinement process continued iteratively until the global error estimator reached a stringent tolerance of $10^{-6}$, ensuring high accuracy in our simulations. For the nonlinear iterations, a convergence tolerance of $10^{-12}$ was applied. Given the non-symmetric nature of the discretized equation for $u$, we employed the \textit{Generalized Minimal Residual Method (GMRES)} to solve the resulting linear systems. To accelerate the convergence of GMRES, \textit{Symmetric Successive Over-Relaxation (SSOR)} was utilized as a preconditioner. In contrast, the linear system of equations originating from the discretization of the $f$-equation was solved using a {direct linear solver}, leveraging its efficiency for this specific part of the problem.

\begin{figure}[H]
    \centering
    \begin{subfigure}[b]{0.48\textwidth} 
        \centering
        \includegraphics[width=\textwidth]{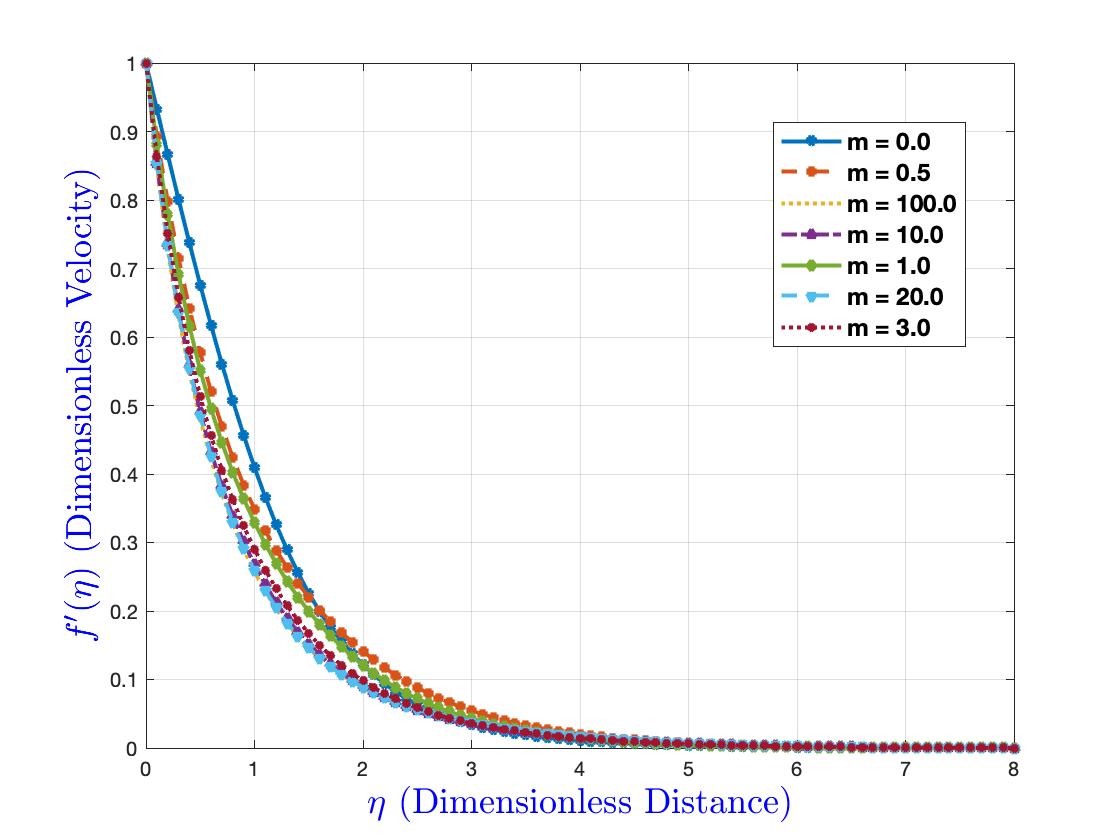}
        \caption{Plot of $f^{\prime}$}
        \label{fig:picard_vnotch_sub}
    \end{subfigure}
    \hfill 
    \begin{subfigure}[b]{0.48\textwidth} 
        \centering
        \includegraphics[width=\textwidth]{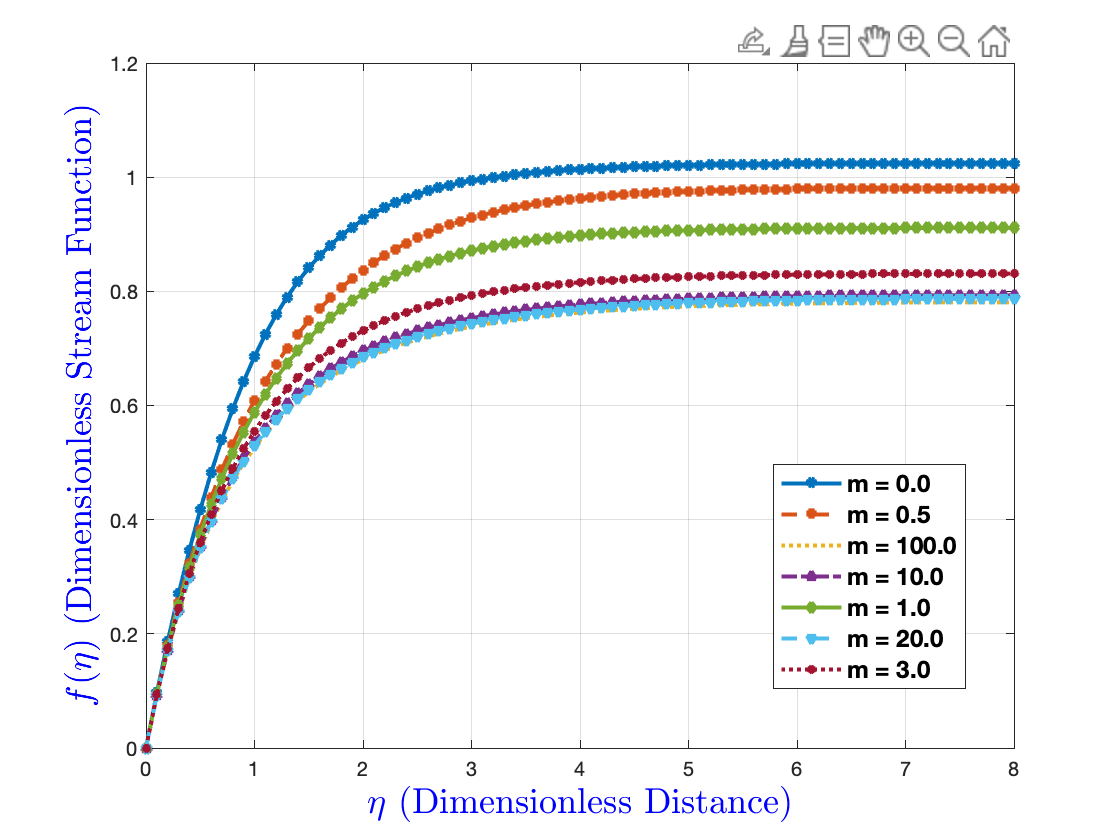}
        \caption{Plot of $f$}
        \label{fig:airy_vnotch_sub}
    \end{subfigure}
    \caption{Plot of $f^{\prime}$ and $f$}
    \label{fig_uprime_f}
\end{figure}

FIgure~\ref{fig_uprime_f} dentoes both dimensioanlless velocity profile ($f^{\prime}$) and dimensionless stream funciton ($f$) for various values of $m$.

\begin{figure}[H]
    \centering
        \centering
        \includegraphics[width=\textwidth]{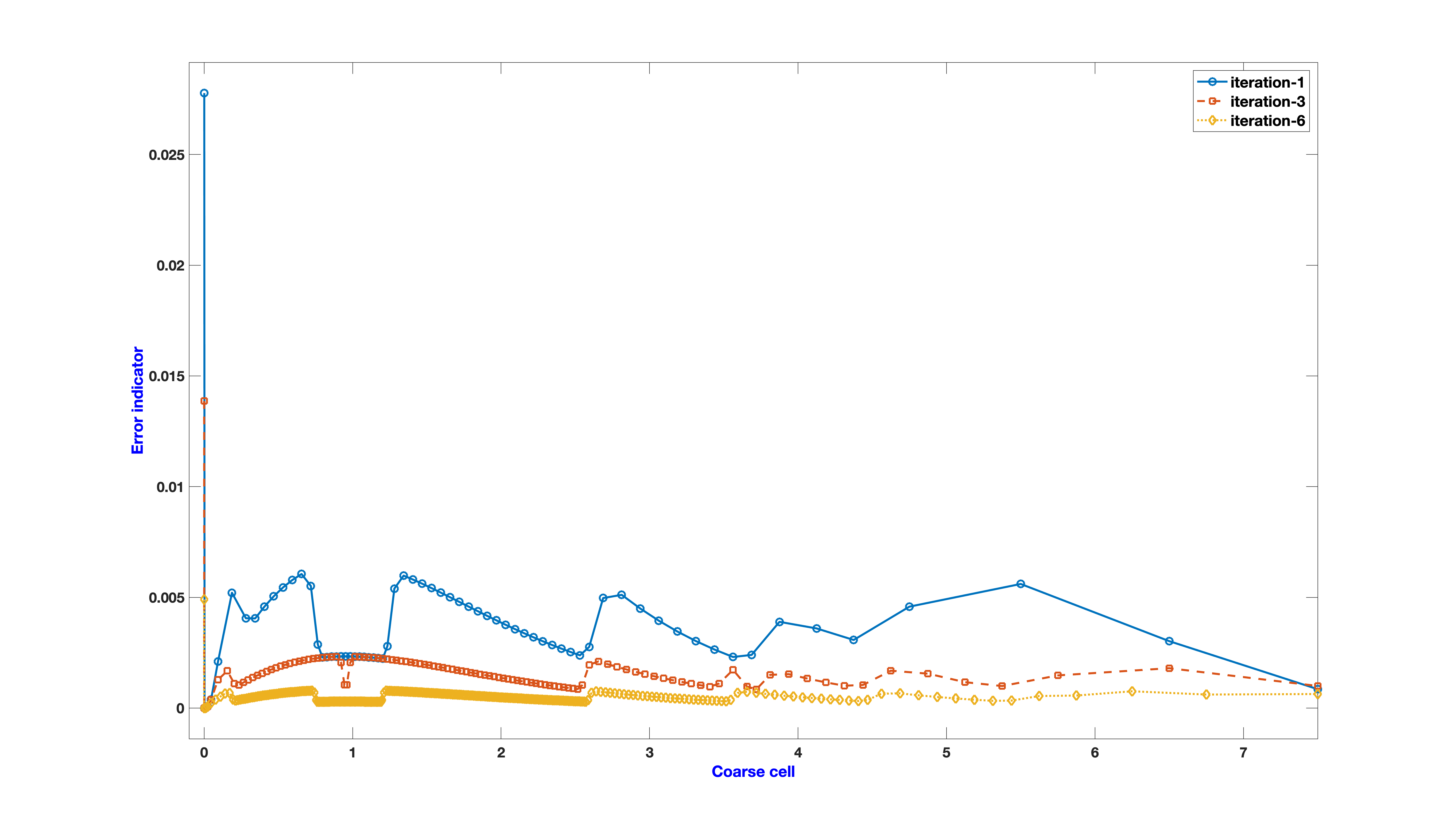}
    \caption{The plot of error indicators across each cell for $m=0$ and $m=1$}
    \label{fig_m1p0}
\end{figure}

Figures~\ref{fig_m1p0} illustrate the {jump in the gradient of $u$ across each computational element} for $m=0$. In this plot, the horizontal axis represents the cell centers, providing the spatial location of each element. The vertical axis indicates the actual magnitude of the gradient jump corresponding to the respective cell centers. It's evident from these figures that the {gradient jumps significantly decreased after the $6^{\text{th}}$ refinement level}. While a substantial reduction was achieved at this stage, our computational process involved {20 refinement steps} to ensure the desired level of accuracy and solution smoothness.

\begin{table}[h!]
    \centering
    \caption{Velocity gradient $f''(0)$ for various values of $m$}
    \label{table_skin_cof}
    \begin{tabular}{cc}
        \toprule
        $m$ & $f''(0)$ \\
        \midrule
        0.0 & 0.67507 \\
        0.2 & 0.93478 \\
        0.5 & 1.1315 \\
        0.8 & 1.2261 \\
        1.0 & 1.2771 \\
        1.5 & 1.3756 \\
        3.0 & 1.4792 \\
        7.0 & 1.5722 \\
        10.0 & 1.604 \\
        20.0 & 1.616 \\
        100.0 & 1.6566 \\
        \bottomrule
    \end{tabular}
\end{table}

Table~\ref{table_skin_cof} presents a comprehensive listing of the {velocity gradient values, specifically $f''(0)$}, for various parameters of $m$. These crucial values were derived directly from the converged solution for $u$, where $u$ is defined as $f'$. To accurately determine the second derivative at $\eta=0$, a {second-order forward difference scheme} was meticulously employed. This approach ensures a precise approximation of the velocity gradient at the boundary, which is essential for our analysis.


\section{Conclusion}\label{conclusions}
This study successfully developed and applied an h-adaptive finite element method to accurately and efficiently solve the highly relevant Falkner-Skan equation, a cornerstone of steady, two-dimensional laminar boundary layer theory. The computational mesh was dynamically refined by integrating the widely recognized Kelly error estimator, which effectively leverages the jump in the gradient across element boundaries for a posteriori error estimation. This adaptive strategy proved instrumental in achieving a precise and computationally economical resolution of the intricate boundary layer behavior inherent to Falkner-Skan flows, which often present significant challenges for non-adaptive numerical techniques due to sharp gradients near the wall.

The robustness and accuracy of the implemented adaptive finite element approach were rigorously demonstrated through its application across a comprehensive range of wedge flow parameters, encompassing both favorable and adverse pressure gradients. A particular emphasis of this investigation was the precise determination of the skin friction coefficient, a critical parameter for quantifying wall shear stress and a direct measure of viscous drag. The obtained numerical solutions consistently provided highly accurate values for this coefficient across all tested conditions, validating the method's capability to capture essential physical quantities. The results presented herein underscore the significant advantages of employing adaptive finite element methods for nonlinear boundary layer problems, showcasing their superior efficiency in allocating computational resources to regions of high interest and their ability to yield reliable solutions where conventional methods might struggle or incur excessive costs. This work thus confirms the adaptive FEM as a powerful tool for advanced boundary layer analysis and opens avenues for its application to more complex fluid dynamic phenomena.

For future work, exploring the application of the Weak Galerkin (WG) finite element method to the Falkner-Skan equation presents a promising avenue, as WG methods offer enhanced flexibility in handling discontinuous approximations and complex geometries, potentially leading to improved accuracy and stability, especially in regions with sharp gradients characteristic of boundary layers  \cite{AAT_SH_SMM2025}.  A compelling future direction involves integrating the {feedforward neural network-based deep learning approximation tools} \cite{martinez2024approximation,venkatachalapathy2023feedforward,venkatachalapathy2023deep,mallikarjunaiah2023deep} with the {adaptive finite element method} presented in this paper. This synergy could lead to the development of a {hybrid deep neural network-driven, mesh-adaptive finite element solver} specifically designed for approximating solutions to complex {nonlinear operators}. Beyond the immediate scope of this study, a promising direction for future research involves the synergistic integration of the adaptive finite element method developed in this paper with the specialized shape functions introduced by Sasikala et al. \cite{sasikala2023efficient}. These advanced shape functions, known for their efficiency in capturing complex solution behaviors, could significantly enhance the performance and accuracy of the adaptive solver. This combined approach holds considerable potential for developing a highly efficient and robust numerical solver specifically tailored for nonlinear fluid flow models, such as those discussed by Yoon et al. \cite{yoon2025steady}. Such an integrated methodology could offer superior capabilities in resolving intricate flow phenomena, improving the fidelity of simulations for a wide range of engineering and scientific applications.

\section*{Acknoledgement}
The author gratefully acknowledges the excellent computing facilities made available by the Department of Mathematics \& Statistics at Texas A\&M University-Corpus Christi, Texas, USA. These resources were crucial in carrying out the computational aspects of this study.

\nocite{*}
\bibliographystyle{plain}
\bibliography{ref_FKEquation}

\end{document}